\documentclass[11pt]{amsart}

\usepackage{amsmath}
\usepackage{amscd}
\usepackage{amsfonts}
\usepackage{graphicx}
\usepackage[usenames,dvipsnames]{color}
\usepackage{color}

\theoremstyle{plain}

\newtheorem{theorem}{Theorem}[section]

\newtheorem{lemma}{Lemma}[section]

\newtheorem{definition}{Definition}[section]

\newtheorem{proposition}{Proposition}[section]

\newtheorem{corollary}{Corollary}[section]

\def \epsilon {\varepsilon}
\def \phi {\varphi}
\def \proof {{\bf Proof}\hspace{.4cm}}

\def \Uu {(U,u)}

\def \AL {A^{\Lambda}}
\def \ALm {A^{\Lambda_m}}
\def \ALa {(A^{\Lambda},a)}

\def \ALam {(A^{\Lambda_m},a_m)}

\def \Kim {K^i_m}

\def \Umum {(U_m,u_m)}
\def \Um {U_m}

\def \gim {\gamma^i_m}

\def \leim {\ell^i_m}

\def \Kim {K^i_m}

\def \Car {Carath\'eodory }

\def \C {\mathbb C}
\def \D {\mathbb D}

\def \R {\mathbb R}

\def \cbar {\overline {\mathbb C}} 
\def \dsharpz {|{\mathrm d}^{\scriptscriptstyle \#}\hspace{-.06cm} z|}

\title{The Carath\'eodory Topology for Multiply Connected Domains I}

\author{Mark Comerford}

\address{Department of Mathematics,
University of Rhode Island,
5 Lippitt Road, Room 102F,
Kingston, RI 02881, USA.
email: {\tt mcomerford@math.uri.edu}}



\begin{document}

\renewcommand{\theenumi}{\emph{\arabic{enumi}.}}
\renewcommand{\labelenumi}{\theenumi}

\maketitle

\begin{abstract}
We consider the convergence of pointed multiply connected domains in the \Car topology. Behaviour in the limit is largely determined by the properties of the simple closed hyperbolic geodesics
which separate components of the complement. Of particular importance are those whose hyperbolic length is as short as possible which we call \emph{meridians} of the domain. We prove a continuity result on convergence of such geodesics for sequences of pointed hyperbolic domains which converge in the \Car topology to another pointed hyperbolic domain. Using this we describe an equivalent condition to \Car convergence which is formulated in terms of Riemann mappings to standard slit domains. 
\end{abstract}

\section{Introduction}
The \Car topology for pointed domains was first introduced in 1952 by \Car \cite{Car} who proved that, for simply connected domains, convergence with respect to this topology is equivalent to convergence of suitably normalized inverse Riemann mappings on compact subsets of the unit disc $\D$. This result is also mentioned by McMullen \cite{McM} who uses it to prove a compactness result for polynomial-like mappings. Our work is also motivated by complex dynamics, in particular the area of non-autonomous iteration where one considers compositions arising from sequences of analytic functions which are allowed to vary. It turns out that in order to prove a non-autonomous version of the classical Sullivan straightening theorem, one must consider the behaviour of multiply connected pointed domains with respect to this topology. 

As we shall see (e.g. in Figure 2 below), one issue is that connectivity is not in general preserved for \Car limits and that some of the complementary components can shrink to a point. This presents problems if one wants to perform quasiconformal surgery on multiply connected domains as certain conformal invariants associated with the domains can become unbounded. One of our ultimate goals, then, will be to find necessary and sufficient conditions for which connectivity is preserved for \Car limits and none of the complementary components of the limit domain is a point (in the finitely connected case, such domains are called \emph{non-degenerate}).

Epstein \cite{Ep} has shown that convergence in the \Car topology is equivalent to convergence of suitably normalized universal covering maps on compact subsets of $\D$ (see Theorem 1.2).
However, it turns out that the limiting behaviour of a sequence of domains of the same connectivity is best understood in terms of certain simple closed hyperbolic geodesics associated with the domain. In Theorem 1.8 we prove the important result 
that if a pointed domain $\Uu$ is a \Car limit of a sequence of pointed domains $\{\Umum\}_{m=1}^\infty$, then every simple closed geodesic of $U$ is a uniform limit of simple closed geodesics of the domains $U_m$ and the corresponding hyperbolic lengths and distances of these geodesics to the basepoints also converge. 

Of particular importance are those geodesics are known as \emph{meridians} which are essentially the shortest simple closed geodesics which separate the complement of the domain in some prescribed way. In Theorem 3.2 we use meridians to prove a version of the above classical result concerning convergence of normalized inverse Riemann mappings for the multiply connected case where we replace the unit disc by suitable slit domains. In the second part of this paper we use meridians to give a solution to our originally stated problem regarding the preservation of connectivity. In fact in Theorem 4.2 we give several equivalent conditions for a family of non-degenerate $n$-connected pointed domains which ensure that any \Car limit is still $n$-connected and non-degenerate.

We begin our exposition with a short resume of the well-known results about the \Car topology.  For the most part we shall be working with the spherical metric $\rm d^{\scriptscriptstyle \#} (\cdot \,, \cdot)$ on $\cbar$ (rather than the Euclidean metric). Recall that the length element for this metric, $\dsharpz$ is given by
\vspace{.25cm}
\[ \dsharpz = \frac{|{\rm d}z|}{1 + |z|^2}\]

and that for an analytic function we have the \emph{spherical derivative}
\vspace{.25cm}
\[ f^\# (z) = \frac{f'(z)}{1 + |f(z)|^2}.\]

\vspace{.25cm}
A {\it pointed domain} is a pair $(U,u)$ consisting of an open  
connected subset $U$ of $\cbar$ (possibly equal to $\cbar$ itself) and a point $u$ in $U$. We say that $\Umum \to \Uu$ in the \Car topology as $m$ tends to infinity if 

\vspace{.2cm}
\begin{enumerate}
\renewcommand{\theenumi}{{\bf \roman{enumi})}}
\renewcommand{\labelenumi}{\theenumi}

\item $u_m \to u$ in the spherical topology;

\vspace{.2cm}
\item For all compact sets $K \subset U$, $K \subset U_m$ for all but finitely many $m$;

\vspace{.2cm}
\item For any {\it connected} (spherically) open set $N$ containing $u$, if $N \subset U_m$ for\\
infinitely many $m$, then $N \subset U$.

\end{enumerate}

We also wish to consider the degenerate case where $U = \{u\}$. In this 
case condition ii) is omitted ($U$ has no interior of which we can take
compact subsets) while condition iii) becomes

\begin{enumerate}

\renewcommand{\theenumi}{{\bf \roman{enumi})}}
\renewcommand{\labelenumi}{\theenumi}
\setcounter{enumi}{2}
\item For any {\it connected} open set $N$ containing $u$, 
$N$ is contained in at most finitely many of the sets $U_m$.

\end{enumerate}

The above definition is a slight modification of that given in the book of McMullen \cite{McM} and much 
of what follows in this section is based on his exposition. However, the original reference for this material goes back to \Car \cite{Car} who in 1952 used an alternative definition which centered on the \emph{\Car kernel} (this approach was also used subsequently by Duren \cite{Dur}). For a sequence of pointed 
domains as above, one first requires that $u_m \to u$ in the spherical topology. If there is no open set 
containing $u$ which is contained in the intersection of all but finitely 
many of the sets $U_m$, one then defines the kernel of the sequence of pointed
domains $\{\Umum\}_{m=1}^\infty$ to be $\{u\}$. Otherwise one then defines the 
\Car kernel as the largest domain $U$ containing $u$ with the 
property ii) above, namely that every compact subset $K$ of $U$ must lie 
in $U_m$ for all but finitely many $m$. It is relatively easy to check that
an arbitrary union of domains with this property will also inherit it. 
Hence a largest such domain does indeed exist. Convergence in this context is 
then defined by requiring that every subsequence of pointed domains has the 
same kernel as the whole sequence. 

It is not too hard to show that this version of \Car convergence is 
equivalent to the first one. In fact, one has the following. 

\begin{theorem}
Let $\{\Umum\}_{m=1}^\infty$ be a sequence of pointed domains and $\Uu$ be another pointed domain where we allow the possibility that $\Uu = (\{u\},u)$. Then the following are 
equivalent:

\vspace{.2cm}
\begin{enumerate}

\item $\Umum \to \Uu$;

\vspace{.2cm}
\item $u_m \to u$ in the spherical topology and $\{\Umum\}_{m=1}^\infty$ has \Car kernel $U$ as does every subsequence; 

\vspace{.2cm}
\item $u_m \to u$ in the spherical topology and, for any subsequence where the complements of the sets $U_m$ 
converge in 
the Hausdorff topology (with respect to the spherical metric), $U$ correspsonds with the connected component of the
complement of the Hausdorff limit which contains $u$ (this component being 
empty in the degenerate case $U = \{u\}$).
 
\end{enumerate}
 
\end{theorem}

It follows easily from the compactness of $\cbar$ combined with the Blaschke selection theorem that, provided we use the spherical rather than the Euclidean metric, any sequence of non-empty closed subsets of $\cbar$ will have a subsequence which converges in the Hausdorff topology. Hence, from above, given 
any family of pointed domains we can always find a sequence in the family which converges in the \Car topology (although the limit pointed domain may well be degenerate). In fact, this convenient fact is the main reason we define things using the spherical topology rather than the more usual Euclidean topology. 

As we remarked earlier, connectivity cannot increase with respect to \Car limits. To be precise, if each $U_m$ above is at most $n$-connected, then so is the limit domain $U$. The reason for this is that 
by {\it 3.} above, complementary components are allowed to merge in the Hausdorff limit, but they cannot split up into more components (see Figure 2 for an illustration of what can happen in this situation). 

Recall that a Riemann surface is called \emph{hyperbolic} if its universal covering space is the unit disc $\D$. From the uniformization theorem, it is well-known that a domain $U \subset \cbar$ is hyperbolic if $\cbar \setminus U$ contains at least three points. For such a domain, the universal covering map allows us to define the \emph{hyperbolic metric} on $U$ which we denote by $\rho_U (\cdot\,, \cdot)$ or just $\rho (\cdot \,, \cdot)$, if the domain involved is clear from the context. Extending this notation slightly, we shall use $\rho_U(z, A)$ or $\rho(z, A)$ to denote the distance in the hyperbolic metric from a point $z \in U$ to a subset $A$ of $U$. Finally, for a curve $\gamma$ in $U$, let us denote the hyperbolic length of $\gamma$ in $U$ by $\ell_U(\gamma)$, or, again when the context is clear, simply by $\ell(\gamma)$.

Often, for the sake of convenience, we shall restrict ourselves to considering domains which are subsets of $\C$ so that the point at infinity is in one of the components of the complement. 
This simplification has the advantage that for a sequence of functions whose ranges lie in domains which are subsets of $\C$ and thus avoid infinity, convergence in the spherical topology is equivalent to the simpler condition of convergence in the Euclidean topology. 

To see why there is little loss of generality in making this assumption, suppose $\Umum$ converges to $\Uu$ with $U$ hyperbolic. Then any Hausdorff limit of the sets $\cbar \setminus U_m$ must contain at least three distinct points since otherwise, $U$ will fail to be hyperbolic. 
We can then apply a M\"obius transformation to $U$ to move these three points to $0$, $1$ and $\infty$. If we now apply the same transformation to the domains $U_m$, then we know that $0$, $1$ and $\infty$ are close to $\cbar \setminus U_m$ for $m$ large. We can therefore choose three points in $\cbar \setminus U_m$ which get moved to $0$, $1$, $\infty$ by a M\"obius transformation which is very close to the identity. It is easy to see from the definition of \Car convergence that this does not affect convergence to the limit domain $\Uu$ and so we have what we want. 

One of the nice features of the \Car topology is that the geometric and 
topological formulations of convergence given above correspond to the 
function-theoretic condition of the local uniform convergence of suitably normalized
covering maps. Of course, in the simply connected case, these are just the inverses of Riemann mappings to the unit disc. We will prove the following result in Section 2.  

\begin{theorem} Let $\{\Umum\}_{m \ge 1}$ be a sequence of pointed hyperbolic domains and for each $m$ let $\pi_m$ be the unique normalized covering map from $\D$ to $U_m$ satisfying $\pi_m(0) = u_m$, $\pi'_m(0) >0$. 

Then $\Umum$ converges in the \Car topology to another pointed hyperbolic domain $\Uu$ if and only if the 
mappings $\pi_m$ converge with respect to the spherical metric uniformly on compact subsets of $\D$ to the covering map $\pi$ from $\D$ to $U$ satisfying $\pi(0)=u$, $\pi'(0) >0$. 

In addition, in the case of convergence, if $D$ is a simply connected subset of $U$ and $v \in D$, then locally defined branches $\omega_m$ of $\pi_m^{\circ -1}$ on $D$ for which $\omega_m(v)$ converges to a point in $\D$ will converge locally uniformly with respect to the spherical metric on $D$ to a uniquely defined branch $\omega$ of $\pi^{\circ -1}$.

Finally, if $\pi_m$ converges with respect to the spherical topology locally uniformly on $\D$ to the constant function $u$, then $\Umum$ converges to $(\{u\},u)$.
\end{theorem}

One of the most important ways to characterize a multiply connected domain is in terms of the simple closed hyperbolic geodesics which separate components of the complement and we will use the tool of homology from complex analysis to classify these curves. The following four results are proved in \cite{Com2}. 

Note that in \cite{Com2} it is always assumed that if a simple closed curve $\gamma$ separates two disjoint closed sets $E$, $F$, then $\infty \in F$. This has the advantage of allowing us to assign a consistent orientation to such a curve so that the winding number $n(\gamma, z)$ is $1$ for all points of $E$ and $0$ for all points of $F$. However, it is obvious that, by applying a suitable M\"obius transformation if needed, we can assume that $E$ and $F$ are any two arbitrary disjoint closed subsets of $\cbar$.  

Another advantage of assuming $\infty \in F$ is that all positively oriented simple closed curves which separate $E$ and $F$ are then in the same homology class and vice versa. If $U \subset \C$, and $\gamma$, $\eta$ are curves in $U$, then we write 
$\gamma \underset{U} \approx \eta$ to denote homology in $U$. 

On the other hand, if we allow $\infty \in U$, then this it is easy to see that there can be curves which separate the complement of $U$ in the same way, but which are not homologous in $U$. This is important for the definition of meridians (see Definition 1.1 below) where we need to take this into account if we wish to consider subdomains of $\cbar$ instead of just subdomains of $\C$.

\begin{theorem}[\cite{Com2} Theorem 2.1] Let $U$ be a domain and suppose we can find disjoint non-empty closed sets $E$, $F$ with $\cbar \setminus U = E \cup F$. Then there exists a piecewise smooth simple closed curve in $U$ which separates $E$ and $F$.
\end{theorem}

For the next three results, we assume the common hypothesis that $U$ is a hyperbolic domain and $E$ and $F$ are closed disjoint non-empty sets neither of which is a point and for which $\cbar \setminus U = E \cup F$. Let us call such a separation of the complement of $U$ \emph{non-trivial}. Also, since we are considering domains which are subsets of $\C$, let us assume that $E$ is bounded and $\infty \in F$. 

\vspace{.15cm}

\begin{theorem}[\cite{Com2} Theorem 2.4]  Let $\tilde \gamma$ be a simple closed curve which separates $E$ and $F$. Then there exists a unique simple closed smooth geodesic $\gamma$ which is the shortest curve in the free homotopy class of $\tilde \gamma$ in $U$ and in particular also separates $E$ and $F$. 

Conversely, given a simple closed smooth hyperbolic geodesic $\gamma$ in $U$, $\gamma$ separates $\cbar \setminus U$ non-trivially and is the unique geodesic in its homotopy class and also the unique curve of shortest possible length in this class. 
\end{theorem}

Note that the fact that $\gamma$ must separate $E$ and $F$ in the first part of the statement follows easily from the Jordan curve theorem and the fact that $\gamma$ is simple and must be homologous in $U$ to $\tilde \gamma$. 
As we will see (e.g. in Figure 1), there may be many geodesics in different homotopy classes which separate $E$ and $F$. However, we can always find one which is as short as possible. 
\vspace{.15cm}

\begin{theorem}[\cite{Com2} Theorem 1.1] Let $U$, $E$ and $F$ be as above. Then there exists a geodesic $\gamma$ which separates $E$ and $F$ and whose length in the hyperbolic metric is as short as possible among all geodesics which separate $E$ and $F$. 
\end{theorem}

Unfortunately, this geodesic need be neither simple nor uniquely defined (see \cite{Com2} for details). However, there does always exist a simple closed geodesic of minimum length which separates $E$ and $F$. 
\vspace{.15cm}

\begin{theorem}[Existence Theorem, \cite{Com2} Theorem 1.3]
There exists a simple closed geodesic $\gamma$ in $U$ which separates $E$ and $F$ and whose hyperbolic length is as short as possible in its homology class and is also as short as possible among all simple closed curves which separate $E$ and $F$. Furthermore, any curve in the homology class of $\gamma$ and which has the same length as $\gamma$ must also be a simple closed geodesic.
\end{theorem}

Note that $\gamma$ is the shortest curve in its homology class which in general includes curves which may not be simple (including possible $\tilde \gamma$ itself). The above statement is a simplified version of the original. In the original, the class of curves which separated $E$ and $F$ by parity was considered and this class is larger than the homology class of $\gamma$ (again, see \cite{Com2} for details). 

Let $\gamma$ be a simple closed smooth hyperbolic geodesic which is topologically non-trivial in $U$, let $\pi : \D \mapsto U$ be a universal covering map and let $G$ be the corresponding group of covering transformations.  Any lift of $\gamma$ to $\D$ is a hyperbolic geodesic in $\D$ and going once around $\gamma$ lifts to a covering transformation $A$ which fixes this geodesic. It is then not hard to see that 
$A$ must be a hyperbolic M\"obius transformation and the invariant geodesic is then $Ax_A$, the axis of $A$. The hyperbolic length of $\gamma$ is then the same as the translation length $\ell(A)$ which is the hyperbolic distance $A$ moves points on $Ax_A$. Note that the quantity $\ell(A)$ does not depend on our choice of lift and is conformally invariant. 

We call a segment of $\eta$ of $Ax_A$ which joins two points $z$, $A(z)$ on $Ax_A$ a \emph{full segment} of $Ax_A$. This discussion and the above result lead to the following definition. 

\begin{definition} Let $U$ be a hyperbolic domain and let $E$, $F$ be any non-trivial separation of $\cbar \setminus U$ as above (where we do not assume that $\infty \notin U$). A simple closed hyperbolic geodesic $\gamma$ in $U$ which separates $E$ and $F$ whose hyperbolic length is as short as possible is called a a \emph{meridian} of $U$ and the hyperbolic length $\ell_U(\gamma)$ is called the \emph{translation length} or simply the \emph{length} of $\gamma$. 
\end{definition}

Note that in Definition 1.5 of \cite{Com2}, a slightly different definition was given where the meridian was defined to be the shortest possible simple closed geodesic in its homology class. As mentioned above, in that paper it was assumed that $\infty \in F$ and in this case the two definitions are equivalent. However, since we wish to consider arbitrary domains in $\cbar$ and not just in $\C$, we need the slightly more general definition above. 

An important special case and indeed the prototype for the above definition is the equator of a conformal annulus and just as the equator is important in determining the geometry of a conformal annulus, meridians are important in determining the geometry of domains of (possibly) higher connectivity. 

\begin{center}
\vspace{.1cm}
\scalebox{0.75}{\includegraphics{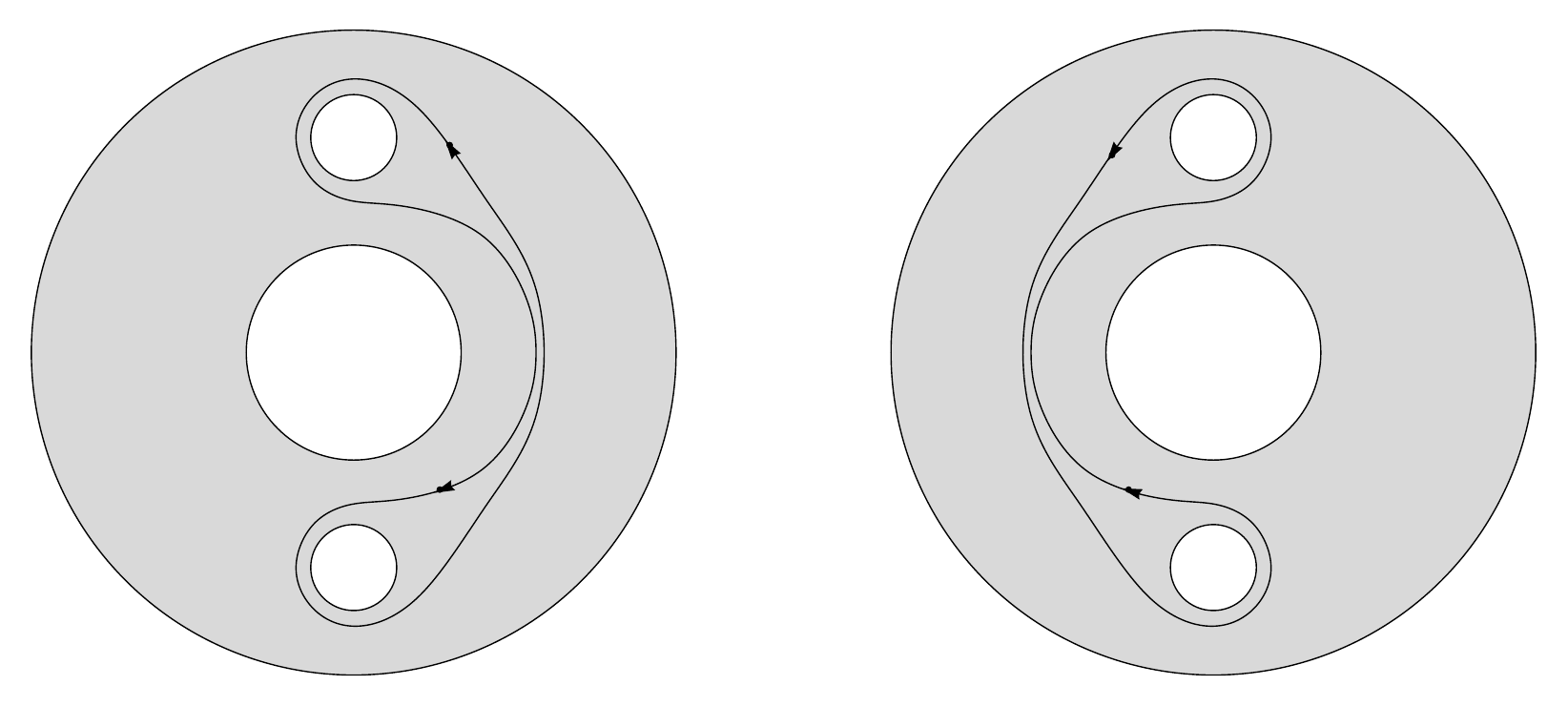}}
\end{center}

The main problem with meridians is that, except in special cases such as an annulus, meridians may not be unique as Figure 1 above shows. The two meridians shown are not homotopic but are in the same homology class and have equal length (see \cite{Com2} Theorem 1.4 for details).

However, if one of the complementary components is connected, then we do have uniqueness. 

\begin{theorem}[\cite{Com2} Theorem 1.5]
If at least one of the sets $E$, $F$ is connected, then there is only one simple closed geodesic $\gamma$ in $U$ which separates $E$ and $F$. In particular, $\gamma$ must be a meridian. In addition, any other geodesic which separates $E$ and $F$ must be longer than $\gamma$.
\end{theorem}

Let us call a meridian as above where at least one of the sets $E$, $F$ is connected a \emph{principal meridian} of $U$. The theorem then tells us that principal meridians are unique. These meridians have other nice properties. For example, they are disjoint and do not meet any other meridians of $U$ (\cite{Com2} Theorem 2.6). 

To see the pathologies which can arise when one takes a limit in the \Car topology, consider Figure 2 below. Note how the connectivity decreases when parts of the complement merge in the limit or are `pinched off' (see Figure 2 below). 

\begin{center}
\scalebox{1.061}{\includegraphics{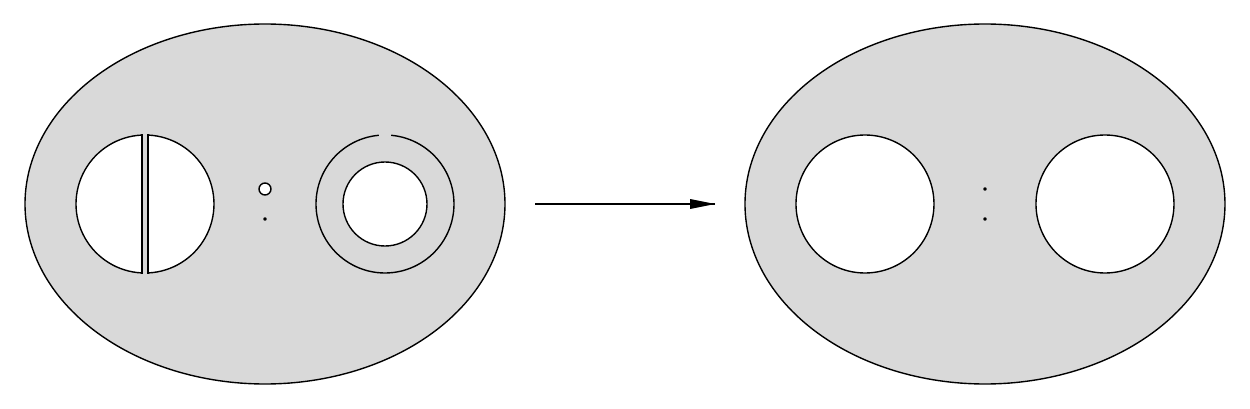}}
\unitlength1cm
\begin{picture}(0.01,0.01)
  \put(-11.4,-.3){\footnotesize$(U_m,u_m)$  }
  \put(-7.42,2.4){\footnotesize $m \to \infty$}
  \put(-3.33,-.3){\footnotesize$(U,u)$  }
  \put(-10.73,1.93){$\scriptstyle u_m$}
  \put(-2.9,1.93){$\scriptstyle u$}
\end{picture}
\end{center}

\vspace{.4cm}

In the above figure the principal meridians which separate the semi-circular shaped complementary components on the left from the rest of $\cbar \setminus U_m$ have lengths which 
must tend to infinity. For the small complementary component in the middle which shrinks to a point, the opposite happens and the principal meridian which separates this component from the rest of the complement has length tending to zero. Finally for the the circular complementary component on the right which is almost swallowed by the circular arc, the principal meridian which separates this component from the rest of the complement will tend to a circle (in fact the equator of a round annulus). However, the distance of this meridian from the base point $u_m$ will tend to infinity. 

The important issue here is that the fact that the limit domain is degenerate and of lower connectivity than the domains of the approximating sequence can be understood entirely in terms of the behaviour of the meridians and in fact of the principal meridians of these domains. Meridians are thus central to understanding the \Car topology in the multiply connected case.

Even though simple closed geodesics can behave badly with respect to limits in the \Car topology, 
we can say something as the theorem below which is one of the main results of this paper shows. Roughly it states that a simple closed geodesic of the limit domain can be approximated by simple closed geodesics of the approximating domains. We say that a sequence of curves $\gamma_m$ converges uniformly to a curve $\gamma$ if we can find parametrizations for all the curves $\gamma_m$ over the same interval which converge uniformly to a parametrization of $\gamma$. 
 
\begin{theorem}
Let $\{\Umum\}_{m=1}^\infty$ be a sequence of multiply connected hyperbolic pointed domains which converges in the \Car topology to a multiply connected hyperbolic pointed domain $\Uu$ (with $U \ne \{u\}$). If $\gamma$ is a simple closed geodesic of $U$ whose length is $\ell$, then we can find simple closed geodesics $\gamma_m$ of each $U_m$ such that if $\ell_m$ is the length of each $\gamma_m$, then:

\begin{enumerate}
\item The hyperbolic distance in $U_m$ from $u_m$ to $\gamma_m$, $d_m = \rho_{U_m}(u_m, \gamma_m)$, converges to $d = \rho_U (u, \gamma)$, the hyperbolic distance in $U$ from $u$ to $\gamma$;

\vspace{.3cm}
\item The simple closed geodesics $\gamma_m$ converge uniformly to $\gamma$ while the corresponding lengths $\ell(\gamma_m)$ converge to $\ell(\gamma)$;

\vspace{.3cm}
\item If $u_m$ lies on $\gamma_m$ for infinitely many $m$, then $u$ lies on $\gamma$.

\end{enumerate}
\end{theorem}

In the case of meridians for a domain, we can say the following. 

\begin{theorem}
Let $\Umum$ and $\Uu$ be as above in Theorem 1.8, 
let $\tilde \gamma$ be a simple closed geodesic of $U$ and let 
$\tilde \gamma_m$ be the curves in each $U_m$ which converge to $\tilde \gamma$ as above. 
For each $m$, let $\gamma_m$ be a meridian of $U_m$ with $\gamma_m \underset{U_m} \approx  \tilde \gamma_m$.  

Then the distances $d_m = \rho_{U_m}(u_m, \gamma_m)$ are uniformly bounded above and the lengths $\ell_m = \ell(\gamma_m)$ are uniformly bounded above and uniformly bounded below away from zero. 
\end{theorem}

\vspace{.1cm}
\begin{theorem}
Again let $\Umum$ and $\Uu$ be as in Theorem 1.8 and suppose $E$, $F$ is a non-trivial separation of $\cbar \setminus U$ into disjoint closed subsets. Then we can find a meridian $\gamma$ which separates $E$ and $F$ and a subsequence $m_k$ and meridians $\gamma_{m_k}$ of $U_{m_k}$ 
such that if $\ell_{m_k}$ is the length of each $\gamma_{m_k}$ and $\ell$ the length of $\gamma$, then:

\begin{enumerate}
\item The hyperbolic distance in $U_{m_k}$ from $u_{m_k}$ to $\gamma_{m_k}$, $d_{m_k} = \rho_{U_{m_k}}(u_{m_k}, \gamma_{m_k})$, converges to $d = \rho_U (u, \gamma)$, the hyperbolic distance in $U$ from $u$ to $\gamma$;

\vspace{.3cm}
\item The meridians $\gamma_{m_k}$ converge uniformly to $\gamma$ while the corresponding lengths $\ell_{m_k}$ converge to $\ell$;

\vspace{.3cm}
\item If $u_{m_k}$ lies on $\gamma_m$ for infinitely many $m$, then $u$ lies on $\gamma$.

\end{enumerate}
Furthermore, if $\gamma$ is a principal meridian of $U$, then {\it 1.}, {\it 2.} and {\it 3.} hold for any subsequence. 

\end{theorem}

An important special case is domains with finite connectivity. We adopt the convention that if $U$ is $n$-connected and $K^1, K^2, \ldots \ldots, K^n$, denote the components of $\cbar \setminus U$, then the last component $K^n$ will always be the unbounded one (note that Ahlfors uses the same convention in \cite{Ahl}). 

For a domain of finite connectivity $n$, one can see using elementary combinatorics that 
there are at most $E(n) := 2^{n-1} - 1$ different ways to separate $\cbar \setminus U$ non-trivially and thus at most this number of meridians which separate the complement of $U$ in distinct ways. One can also show that there are at most $P(n) := \min\{n,E(n)\}$ principal meridians. If we can find $P(n)$ principal meridians, let us call such a collection the \emph{principal system of meridians} or simply the \emph{principal system} for $U$. If we can find a full collection of $E(n)$ meridians, let us call such a collection an \emph{extended system of meridians} or simply an \emph{extended system} for $U$. 

If $n \le 3$, than any meridians of $U$ which exist must be principal. The first case where we can have meridians which are not principal is $n=4$ as we see in Figure 1. Finally, as the principal meridians are always disjoint and in different and non-trivial homotopy classes, they form a geodesic multicurve in the sense of Definition 3.6.1 of \cite{Hub}. However, except when $n = 2$ or $3$, this multicurve will not be maximal (see Theorem 3.6.2 of \cite{Hub}). On the other hand, the meridians of an extended system may well intersect and so will not in general be a multicurve at all. 

A finitely connected domain $U$ is called \emph{non-degenerate} if none of the components of $\cbar \setminus U$ is a point. The principal meridians are precisely those meridians which can fail to exist if some of the complementary components are points and it is not hard to show the following.  

\begin{proposition}[\cite{Com2} Proposition 3.1]
If $U$ is a domain of finite connectivity $n \ge 2$, then $U$ has at least $E(n) - P(n)$ meridians and any principal meridians of $U$ which exist are uniquely defined. Furthermore, the following are equivalent:

\begin{enumerate}
\item $U$ is non-degenerate;

\vspace{.3cm}
\item $U$ has $P(n)$ principal meridians;

\vspace{.3cm}
\item $U$ has $E(n)$ meridians in distinct homology classes.

\end{enumerate}

\end{proposition}

If $U$ is a non-degenerate $n$-connected domain and $\Gamma = \{\gamma^i, 1 \le i \le E(n)\}$ is an extended system for $U$, we shall adopt the convention that the first 
$P(n)$ meridians are always those of the principal system and that for $1 \le i \le P(n)$, $\gamma^i$ separates $K^i$ from the rest of $\cbar \setminus U$. Let us denote the lengths of the meridians of $\Gamma$ by $\ell^i$, $1 \le i \le E(n)$. For a pointed domain $\Uu$, we will also need consider the distances $d^i := \rho(u, \gamma^i)$, $1 \le i \le E(n)$ from the base points to these meridians. 

The collection of numbers $\ell^i$ and $d^i$, $1 \le i \le E(n)$, we shall refer to as the \emph{lengths} and \emph{distances} of $\Gamma$ respectively and naturally we can make similar definitions for a principal system. Note that the lengths are independent of the choice of meridians for the system, but, except for the principal meridians, the distances in general are not. However, this will not be too much of a problem as we see from Theorem 1.9.

Let us call a meridian \emph{significant} if it is a limit of meridians for a subsequence as above. If the domain is finitely connected and non-degenerate, let us call a system of meridians a \emph{significant system} if each meridian in the system is significant. We have the following useful corollary.

\begin{corollary}
Let $n \ge 2$ and let $\{\Umum\}_{m=1}^\infty$ be a sequence of $n$-connected hyperbolic pointed domains which 
converges in the \Car topology to a non-degenerate hyperbolic $n$-connected pointed domain $\Uu$ with $U \ne \{u\}$. Then we can find a significant system of meridians for $\Uu$. 

Furthermore, if all the domains $U_m$ are also $n$-connected and non-degenerate and for each $m$ we let $\Gamma_m = \{\gim, 1 \le i \le E(n)\}$ be any system of meridians for $U_m$, then the distances $d^i_m = \rho(u_m , \gim)$ are bounded above while the lengths $\leim = \ell_{U_m}(\gim)$ are bounded above and below away from zero.  These bounds are uniform in $m$ and independent of our choice of the systems $\Gamma_m$.  



\end{corollary}

Theorems 1.2, 1.8, 1.9 and 1.10 and Corollary 1.1 will be proved in Section 2. In Section 3 we will present some applications including a version of Theorem 1.2 stated in terms of Riemann mappings to slit domains instead of universal covering maps.

\section{Convergence of Geodesics and Meridans}

Starting with Theorem 1.2, we prove the theorems stated in the previous section, together with some supporting results. For a family of M\"obius transformations $\Phi = \{\phi_\alpha, \alpha \in A\}$, we say that $\Phi$ is \emph{bi-equicontinuous on $\cbar$} if both $\Phi$ and the family $\Phi^{\circ -1} =  \{\phi_\alpha^{\circ -1}, \alpha \in A\}$ of inverse mappings are uniformly Lipschitz families on $\cbar$ (with respect to the spherical metric).

{\bf Proof of Theorem 1.2 \hspace{.4cm}} A proof of most of this result can be found in the Ph.D. thesis of Adam Epstein \cite{Ep} and the proof is similar to the better known special case where all the domains involved are discs and the mappings $\pi_m$ are then Riemann maps. A proof of the disc case can be found in {\Car}'s original exposition \cite{Car}. 

In order to extend Epstein's results to a full proof, we need to show in the non-degenerate case that if $\Umum$ converges to $\Uu$ with $U$ hyperbolic, then the covering maps $\pi_m$ give a normal family on $\D$ and that any limit function must be non-constant. Note that in the non-degenerate case, we may (if we like) assume that $U \subset \C$ so that the sequence $u_m$ is bounded in the case of either \Car convergence or convergence of normalized covering maps and so convergence in the spherical topology is equivalent to convergence in the Euclidean topology. Lastly, in the degenerate case we need to show that $\Umum$ converges to $(\{u\},u)$ as stated. 

Dealing first with the non-degenerate case, 
since $U$ is hyperbolic, it then follows from the Hausdorff version of \Car convergence that we can find $\delta >0$ such that for every $m$ large enough $\cbar \setminus U_m$ contains at least three points which are at least distance $\delta$ away from each other in terms of the spherical metric. The reason for this is that if this were not true we could find a subsequence which converged to a domain which was $\cbar$ with one or two points removed, both of which are impossible (note that this argument also shows that if $\Umum \to \Uu$ with $U$ hyperbolic, then $U_m$ must be hyperbolic for $m$ large enough). 

Using Theorem 2.3.3 on page 34 of \cite{Bear}, we can postcompose by a bi-equicontinuous family of M\"obius transformations and apply Montel's theorem to conclude that the covering maps $\pi_m$ give a normal family (in the spherical topology) on $\D$. 
Since $U \ne \{u\}$, it follows from i) and ii) of \Car convergence and applying the Koebe one-quarter theorem to branches of inverse maps on a suitable disc about $u$ in $U$ that all limit functions must be non-constant and this completes the proof in the non-degenerate case.

For the degenerate case, suppose $\pi_m$ converges to the constant function $u$ locally uniformly on $\D$ but $\Umum$ does not converge to $(\{u\},u)$. By Theorem 1.1, using the Hausdorff version of \Car convergence, we can find a subsequence $(U_{m_k}, u_{m_k})$ so that these pointed domains converge to a pointed domain $(\tilde U, u)$ where $\tilde U$ is open with $u \in \tilde U$. From the non-degenerate case, it would then follow that $\pi_{m_k}$ would converge locally uniformly on $\D$ to $\tilde \pi$, the normalized covering map for $V$ and with this contradiction the proof is complete. $\Box$

As one might suspect from the statement of Theorem 1.2, it does not follow that if $\Umum$ converges to a degenerate pointed domain $(\{u\},u)$, then the normalized covering maps $\pi_m$ as above must converge locally uniformly to $u$ on $\D$. The basic reason this fails is that it is possible that the sequence $\{\pi_m\}_{m=1}^\infty$ does not give a normal family and we now give a counterexample which exhibits this behaviour. 

For each $m \ge 1$, let $U_m = {\mathrm A}(0,\tfrac{1}{m^3}, m)$, $u_m= \tfrac{1}{m}$ and consider the sequence of pointed domains $(U_m, u_m)$. This sequence clearly tends to $(\{0\},0)$ and if the family of covering maps had a locally convergent subsequence $\pi_{m_k}$, then it would follow from Rouch\'e's theorem and local compactness as argued by Epstein that  $\pi_{m_k}$ must tend to the constant function $0$ locally uniformly on $\D$. However, it is easy to see that the annulus ${\mathrm A}(0,\tfrac{1}{m^2}, 1)$ has uniformly bounded hyperbolic diameter in $U_m$ (as it has the same equator and half the modulus of the larger annulus) and thus in $U_m$ by the Schwarz lemma for the hyperbolic metric. Since this annulus contains the base point $u_m =  \tfrac{1}{m}$ (which actually lies on its equator), it follows that we can find points $z_m$ within bounded hyperbolic distance of $0$ in $\D$ with $\pi_m(0) = 1$. With this contradiction, we see that the sequence of covering maps cannot have a convergent subsequence and in particular cannot converge as required. 

As McMullen {\cite{McM} Page 67 Theorem 5.3) remarks in the disc case, we can move the base points for a convergent sequence of pointed discs by a uniformly bounded hyperbolic distance without affecting whether or not the sequence converges.  The proof of this fact is a straightforward application of Theorem 1.2.

\begin{corollary} Let $\{\Umum\}_{m=1}^\infty$ be a sequence of pointed hyperbolic domains which converges to $\Uu$ with $U$ hyperbolic. 

\begin{enumerate}

\item If $w_m \in U_m$ for each $m$, $w \in U$ and  $w_m \to w$ as $m \to \infty$ , then $(U_m, w_m)$ converges to $(U,w)$.

\vspace{.2cm}
\item If $w_m \in U_m$ for each $m$ and we can find $d >0$ independent of $m$ so that $\rho_{U_m}(u_m, w_m) \le d$, then we can find $w \in U$ and a subsequence of the sequence of pointed domains $\{(U_m, w_m)\}_{m=1}^\infty$ which converges to $(U,w)$. 

\end{enumerate}

\end{corollary}

The following Lemma will be very useful to us in proving a number of results, especially Theorem 1.8. For $z \in \cbar$ and $r > 0$, let us denote the open spherical disc of radius $r$ about $z$ by ${\mathrm D^\#}(z,r)$.

\begin{lemma}Suppose $\{\Umum\}_{m=1}^\infty$ is a sequence of pointed multiply connected domains which converges to a pointed multiply connected domain $\Uu$ in the \Car topology (where we include the degenerate case $U = \{u\}$) and suppose in addition that the complements $\cbar \setminus \Um$ converge in the Hausdorff topology (with respect to the spherical metric on $\cbar$) to a set $K$.
Then $\partial U \subset K$.
\end{lemma}

\proof Suppose first that we are in the degenerate case where $U = \{u\}$. By iii) of \Car convergence in the degenerate case, for any $0 < \epsilon \le \pi$, ${\mathrm D^\#}(u, \epsilon)$ contains points of $\cbar \setminus U_m$ for all but finitely $m$ and on letting $\epsilon \to 0$, we see that $\partial U = \{u\} \subset K$ as desired.  

Now suppose that $U \ne \{u\}$. By the Hausdorff version of \Car convergence, $U$ is the component of the complement of this Hausdorff limit which contains the point $u$. So suppose the conclusion fails. If $z$ were a point in $\partial U$ which missed $K$, we could find a spherical disc ${\mathrm D^\#}(z, \delta)$ of some radius $\delta >0$ about $0$ which missed $K^m$ for all $m$ sufficiently large. Then $U \cup {\mathrm D^\#}(z, \delta)$ would be a connected set which missed $K$ and since this set contains $z \notin U$, this would contradict the maximality of $U$ as a connected component of $\cbar \setminus K$ whence $\partial U \subset K$ as desired. $\Box$

\begin{center}
\scalebox{1.07}{\includegraphics{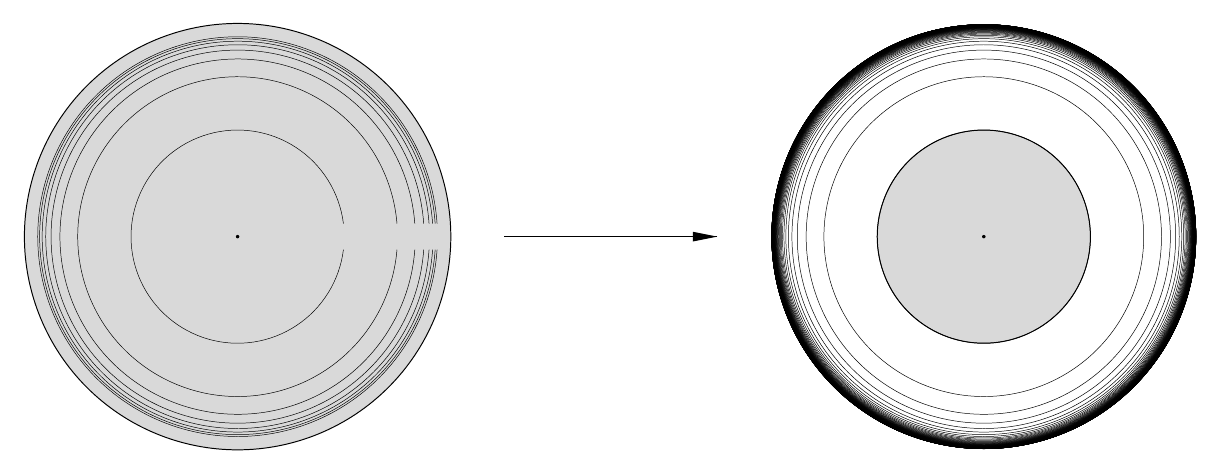}}
\unitlength1cm
\begin{picture}(0.01,0.01)
  \put(-11.4,-.3){\footnotesize$(U_m,0)$  }
  \put(-7.4,2.8){\footnotesize$m \to \infty$}
  \put(-3.1,-.3){\footnotesize$(U,0)$  }
\end{picture}
\end{center}

\vspace{.7cm}

The reader might wonder if, in the case where the limit domain $U$ above was $n$-connected, would this force the set $K$ to have exactly $n$ components. However this is false as the following counterexample depicted in Figure 2 shows. For each $m$ and each $ 2 \le i \le m$, let $\Kim$ be the circle ${\mathrm C}(0,1 - 1/i)$ with an arc of height $1/m$ centered about $1-1/i$ removed. If we then set $U_m = \D \setminus \bigcup_{2 \le i \le m} \Kim$, then the pointed domains $(U_m, 0)$ converge to $({\mathrm D}(0,1/2), 0)$ while their complements converge to the set $\bigcup_{i>=2}{\mathrm C}(0,1-1/i) \bigcup (\cbar \setminus \D)$ which clearly has infinitely many components. 
\vspace{.2cm}

{\bf Proof of Theorem 1.8} \hspace{.4cm} Suppose that $\Umum$ converges to $\Uu$ as in the statement in which case we know from Theorem 1.2 that the corresponding normalized covering maps $\pi_m$ converge uniformly on compact sets of $\D$ to the normalized covering map $\pi$ for $U$. 

Let $\eta$ be a lifting of $\gamma$ to $\D$ which is as close as possible to $0$ and let $A$ be the corresponding hyperbolic covering transformation whose axis is $\eta$. Let $\sigma = [a,b]$ be a full segment of $\eta$ with $b = A(a)$ and which contains the closest point on $\eta$ to $0$. 

Let $\epsilon > 0$ be small and take a small hyperbolic disc $D$ in $\D$ about $a$ of radius $\epsilon$. $\tilde D = A(D)$ is then a disc of hyperbolic radius $\epsilon$ about $b= A(a)$ and we have that $\pi \equiv \pi \circ A$ on $D$. If we let 
$R$ be an elliptic rotation of \emph{angle} $\pi$ about $b$, then $R$ and hence $R \circ A$ cannot belong to the group of covering transformations of $U$ as this would violate the local injectivity at $b$ of the covering map $\pi$. 

From above, the difference $\pi - \pi \circ R \circ A$ cannot be identically zero on $\D$ as otherwise it would follow from the monodromy theorem that $R \circ A$ and hence $R$ would belong to the group of covering transformations. $\pi - \pi \circ R \circ A$ then has an isolated zero at $a$ and so, for $\epsilon$ small enough, 
$\pi - \pi \circ R \circ A$ is non-zero on the boundary of $D$. By the local uniform convergence of $\pi_m$ to $\pi$ on $\D$, if we apply Rouch\'e's theorem to $D$, we see that for $m$ large enough, we can find points $a_m \in D$ and $b_m = R(A(a_m)) \in R(A(D))$ with $\pi_m(a_m) = \pi_m(b_m)$. Since $\epsilon$ was arbitrary, if we let $\sigma_m'$ be the geodesic segment between $a_m$ and $b_m$, then $a_m \to a$, $b_m \to b$ and $\sigma_m' \to \sigma$ as $m \to \infty$.

Thus we can find covering transformations $A_m$ of $\D$ for $\pi_m$ with $A_m(a_m) = b_m$. Now let $\gamma_m'$ be the image of $\sigma_m'$ under $\pi_m$. Note that since $\gamma$ is simple while $\sigma_m'$ is very close to $\sigma$, it follows again from the convergence of $\pi_m$ to $\pi$ that, moving $a_m$ and $b_m$ slightly closer together along $\sigma_m$ if needed by an amount which will tend to $0$ as $m \to \infty$, we can assume that there are no points of self-intersection on $\gamma_m'$. $\gamma_m'$ is then a simple closed curve which is a geodesic except at possibly one point where it is not smooth (i.e. there may be a corner). 

$\gamma$ is also a simple closed curve and as before we will let $E$ and $F$ denote the intersection of $\cbar \setminus U$ with each of the two complementary components of $\gamma$ and assume that $E$ is bounded and $\infty \in F$. 
Since $\gamma$ is a geodesic, each of $E$ and $F$ must contain at least two points in view of the second part of Theorem 1.4. If we let $z, w \in E$ be two such points, 
then we may assume that they are in $\partial E \subset \partial U$.
As $\gamma_m'$ is very close to $\gamma$, the winding number of $\gamma_m'$ about $z$ will be close to that of $\gamma$ about $z$ and the same will be true for $w$. As the curves $\gamma_m'$ are simple, $z$ and $w$ are then inside $\gamma_m'$ for $m$ large and it then follows from Lemma 2.1 that for $m$ large enough there are at least two points of $\cbar \setminus U_m$ inside $\gamma_m'$ while the same argument shows that we may also assume the same about the outside of $\gamma_m'$. 

$\gamma_m'$ is thus a simple closed curve which separates $\cbar \setminus U_m$ non-trivially and we may now apply Theorem 1.4 for $m$ large enough to find a simple closed geodesic $\gamma_m$ which is homotopic in $U_m$ to $\gamma_m'$.  By lifting the homotopy, we can then find a lifting of $\gamma_m$ which coincides with the axis of $A_m$ which we will denote by $\eta_m$.

\vspace{.2cm}
\begin{center}
\scalebox{1.037}{\includegraphics{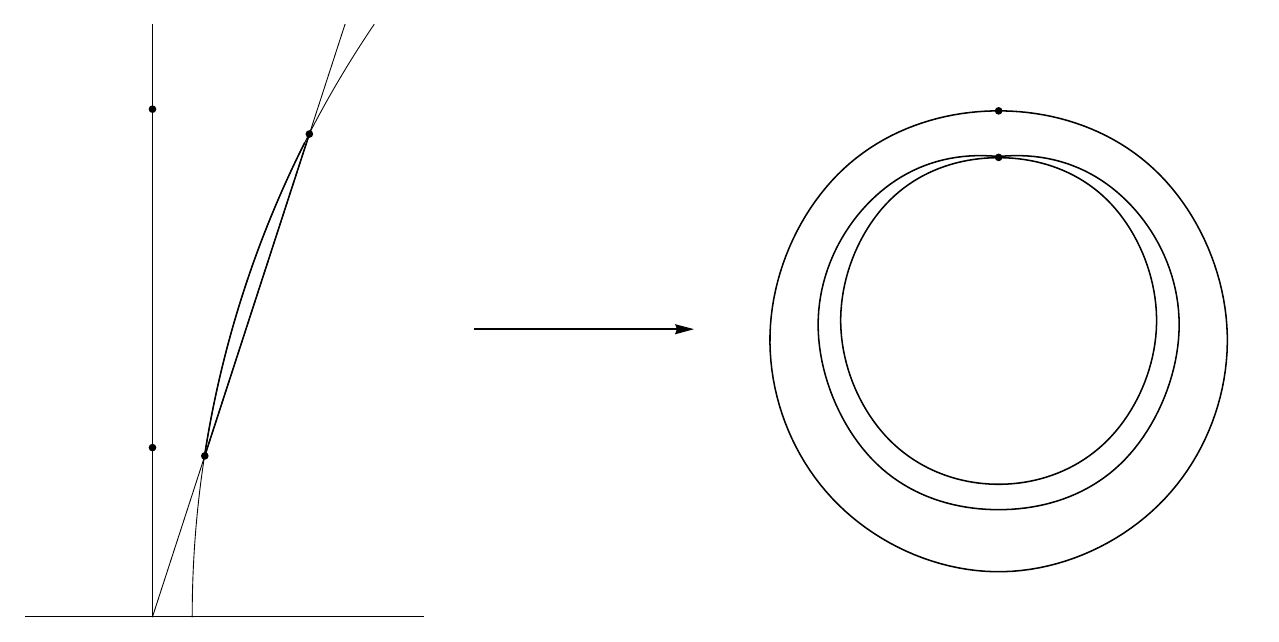}}
\unitlength1cm
\begin{picture}(0.01,0.01)
   \put(-12.2,6.7){\tiny$\eta_m$}
   \put(-12.45,3.7){\tiny$\sigma_m$ }
   \put(-11.55,3.6){\tiny$\sigma_m'$}
   \put(-10.85,3.4){\tiny$\tau_m$  }
   \put(-12.45,2){\tiny$s_m$}
   \put(-12.45,5.5){\tiny$t_m$}
   \put(-11.87,1.9){\tiny$a_m$  }
   \put(-10.8,5.3){\tiny$b_m$  }
   \put(-7.66,3.5){$\pi_m$}
   \put(-3.22,.5){\tiny$\gamma_m$}
   \put(-3.22,1.12){\tiny$\gamma_m'$}
   \put(-3.52,1.9){\tiny$\pi_m(\tau_m)$}
   \put(-3.22,4.85){\tiny$z_m$}
   \put(-3.52,5.85){\tiny$\pi_m(s_m)$}
\end{picture}
\end{center}

\vspace{0cm}

The circle which passes through $a_m$, $b_m$ and the fixed points of $A_m$ is invariant under $A_m$ and its image under $\pi_m$ is a smooth closed curve which in particular has no corner at the point $\pi_m(a_m) = \pi_m(b_m)$ which we will call $z_m$ (this is easiest to see in the upper half-plane model as in the figure above where we let $0$ and $\infty$ be the fixed points of $A_m$ and the imaginary axis the axis of $A_m$ where this circle corresponds to a ray connecting $0$ to $\infty$). Let $\tau_m$ be the segment of this circle which passes through $a_m$ and $b_m$.

Now $\sigma_m'$ is very close to $\sigma$ for $m$ large and since $\pi_m$ converges locally uniformly on $\D$ to $\pi$, the derivatives $\pi_m'$ converge locally uniformly to $\pi'$. Thus 
the difference between the angles of the two tangents to $\gamma_m'$ at the corner at $z_m$ will be very small and will tend to $0$ as $m$ tends to infinity. Since the covering maps $\pi_m$ are angle-preserving, we can say the same about the angles of the tangents at the two endpoints $a_m$, $b_m$ of $\sigma_m'$ (note that this applies in the upper half-plane picture rather than that for the unit disc). 

Now the image $\pi_m(\tau_m)$ of the above invariant circle under $\pi_m$ is a smooth curve and it is clear from the picture above that $\sigma_m'$ must lie on one side of $\tau_m$. It then follows from above that $\sigma_m'$ must then be very close to $\tau_m$. However, since the hyperbolic distance between $a_m$ and $b_m$ is bounded below, this can only happen if $\sigma_m'$ is very close to a segment $\sigma_m$ of the the axis $\eta_m$ of $A_m$ which connects points $s_m$, $t_m$ with $t_m = A_m(s_m)$ (again this is easiest to see in the upper half-plane picture above).  

Hence $\sigma_m$ is very close to $\sigma'_m$ which in turn is very close to $\sigma$ and since all three of these are geodesic segments, their lengths in the hyperbolic metric of $\D$ will also be close. Since these segments are all mapped to simple closed curves by their corresponding covering maps, this gives us \emph{1.} and the second part of \emph{2.} immediately while the rest of \emph{2.} follows on applying the local uniform convergence of $\pi_m$ to $\pi$. 
Finally, \emph{3.} follows immediately from \emph{2.}, which finishes the proof. $\Box$

We remark that the proof above relied mostly on the convergence of normalized covering maps. The only place where we needed \Car convergence directly was for Lemma 2.1 which was used just once to show that the curve $\gamma_m'$ separated $\cbar \setminus U_m$ non trivially. We turn now to proving Theorem 1.9. We first need a lemma from \cite{Com2}. Note that the original version of this lemma was for subdomains of $\C$ where two positively oriented curves separate the complement of $U$ in the same way if and only if they are homologous in $U$. As, usual, however, any hyperbolic domain in $\cbar$ can be mapped to a hyperbolic domain in $\C$ using a M\"obius transformation. 

\begin{lemma}[\cite{Com2} Lemma 2.3] Let $U \subset \cbar$ be a hyperbolic domain and let $\gamma_1$, $\gamma_2$ be two simple closed geodesics in $U$ which separate $\cbar \setminus U$ in the same way and suppose that one of these curves lies in the closure of one of the complementary components of the other. Then $\gamma_1 = \gamma_2$. 
\end{lemma}

{\bf Proof of Theorem 1.9 \hspace{.4cm}} Let $E$ and $F$ be the subsets of $\cbar \setminus U$ separated by $\tilde \gamma$ and, as usual, we assume that $\infty \in F$. Now let $\tilde \gamma_m$ be the geodesics in $U_m$ which converge to $\tilde \gamma$ as in Theorem 1.8. Now for each $m$, let $\gamma_m$ be a meridian which separates the complement $\cbar \setminus U_m$ in the same way as $\tilde \gamma_m$. Applying a M\"obius transformation if needed to map one of the points of $\cbar \setminus U_m$ to $\infty$, we can conclude by Lemma 2.2 above that $\gamma_m$ must meet $\tilde \gamma_m$ and it follows from Theorem 1.8 that the hyperbolic distances $\rho_{U_m}(u_m, \gamma_m)$ must be uniformly bounded above. 

By Theorems 1.6 and 1.8, the lengths $\ell_m$ are obviously bounded above. To see that they must be bounded below, for each $m$ let $\pi_m$ be the normalized universal covering map for $U_m$ and let $\pi$ be the normalized universal covering map for $U$. By Theorem 1.2, $\pi_m$ then converges locally uniformly on $\D$ to $\pi$. 

Now for each $m$, let $\sigma_m$ be a full segment of a lift of $\gamma_m$ which is as close as possible to $0$. The segments $\sigma_m$ are all within bounded distance of $0$ and have uniformly bounded hyperbolic lengths. It then follows that the lengths of these segments and hence the curves $\gamma_m$ must be bounded below away from $0$ since otherwise, by the local uniform convergence of $\pi_m$ to $\pi$ we would obtain a contradiction to the fact that $\pi$ as a covering map must be locally injective. $\Box$

\vspace{.2cm}
{\bf Proof of Theorem 1.10 \hspace{.4cm}} Let $\tilde \gamma$ be a meridian in $U$ which separates $E$ and $F$ which exists by virtue Theorem 1.6. By the discussion at the end of page 3 and the start of page 4 about post-composing with suitably chosen M\"obius transformations, we can assume that $\infty \in F$ and also that $\infty \notin U_m$ for every $m$ which allows us to make use of homology as a tool to completely describe how a simple closed curve separates the complements of these domains. 

By Theorem 1.8, we can find a sequence of geodesics $\tilde \gamma_m$ which tends to $\tilde \gamma$. By Theorem 1.6 again, we can then find meridians $\gamma_m$ in the homology class of each $\tilde \gamma_m$. By Theorem 1.9, the associated distances $d_m$ for the curves $\gamma_m$ are uniformly bounded above while the lengths $\ell_m$ are again uniformly bounded above and bounded below away from zero.

If we now let $\pi_m$ and $\pi$ be the normalized covering maps for each $U_m$ and $U$ respectively, then again $\pi_m$ converges locally to $\pi$ by Theorem 1.2. As above, we can then find full segments $\sigma_m$ of liftings of each $\gamma_m$ which are a uniformly bounded hyperbolic distance from $0$ which are the axes of hyperbolic M\"obius transformations $A_m$ of bounded translation length. It then follows  that we can find a subsequence $m_k$ for which the corresponding segments $\sigma_{m_k}$ converge to a geodesic segment $\sigma$ which must have positive length otherwise we again obtain a contradiction to the local injectivity of $\pi$ as at the end of the proof of Theorem 1.9. If we then set $\gamma = \pi(\sigma)$, then $\gamma$ is a closed hyperbolic geodesic of $U$ and the meridians $\gamma_{m_k}$ must converge to $\gamma$. 

We still need show that this geodesic is simple and a meridian which separates $\cbar \setminus U$ into the same sets $E$, $F$ as 
$\tilde \gamma$ does. As the curves $\tilde \gamma_{m_k}$ , $\gamma_{m_k}$ converge to $\tilde \gamma$, $\gamma$ respectively,  by ii) of \Car convergence, $\tilde \gamma_{m_k}$ , $\gamma_{m_k}$ are bounded away from the boundaries $\partial U_{m_k}$. Thus if $z \in \partial U$ then, by Lemma 2.1, for $k$ large we can find a point $z_{m_k}$ which is very close to $z$ and thus in the same complementary region of $\gamma_{m_k}$. The same argument allows us to make a similar conclusion for the curves $\tilde \gamma_{m_k}$. Hence for $z \in \partial U$ and $k$ large, by homology in $U_{m_k}$, 
\vspace{.3cm}
\[n(\gamma_{m_k}, z) = n(\gamma_{m_k}, z_{m_k}) = n(\tilde \gamma_{m_k}, z_{m_k}) = n(\tilde \gamma_{m_k}, z).\]
\vspace{-.3cm}

It then follows from the above that for $k$ large enough $\tilde \gamma_{m_k}$ and $\gamma_{m_k}$ are homologous in $U$ as well as in $U_{m_k}$ (it is not hard to see that there is sufficient generality in considering winding numbers around points only of $\partial U$ rather than all of $\cbar \setminus U$). Also, by the uniform convergence of the curves $\tilde \gamma_{m_k}$ and $\gamma_{m_k}$ to $\tilde \gamma$ and $\gamma$ respectively, these curves eventually lie in $U$ and are homologous in $U$ to $\tilde \gamma$ and $\gamma$ respectively. Hence, for large $k$ we have 
\vspace{.2cm}
\[ \gamma \underset{U} \approx \gamma_{m_k} \underset{U} \approx \tilde \gamma_{m_k} \underset{U} \approx \tilde \gamma.\]

\vspace{0cm}
Thus $\gamma$ is homologous in $U$ to $\tilde \gamma$ and since $\tilde \gamma$ is a meridian, the length of $\gamma$ cannot be smaller than that of $\tilde \gamma$. On the other hand, by the convergence of the curves $\gamma_{m_k}$ to $\gamma$ using universal covering maps above, and the fact that the curves $\gamma_{m_k}$ are meridians, it follows 
that $\gamma$ cannot be longer than $\tilde \gamma$ either. By Theorem 1.6, $\gamma$ is then a meridian which separates $E$ and $F$ and in particular simple which completes the proof. $\Box$

{\bf Proof of Corollary 1.1} \hspace{.4cm} The existence of a significant system of meridians for $\Uu$ is immediate in view of Theorem 1.10. Now let $\gamma^i$, $1 \le i \le E(n)$ be any extended system of meridians for $U$ and let $\gim$ be the curves which converge to each $\gamma^i$ as in Theorem 1.8. By Lemma 2.1, we see that for $m$ large enough, the curves $\gim$ give different separations of the complement $\cbar \setminus U_m$.  This implies that for $m$ large enough, any meridian of $U_m$ separates the complement of $U_m$ in the same way as one of the curves $\gim$ and the uniform bounds on the distances and lengths of the system then follows from Theorem 1.9. $\Box$

\vspace{.2cm}
\begin{center}
\scalebox{1.04}{\includegraphics{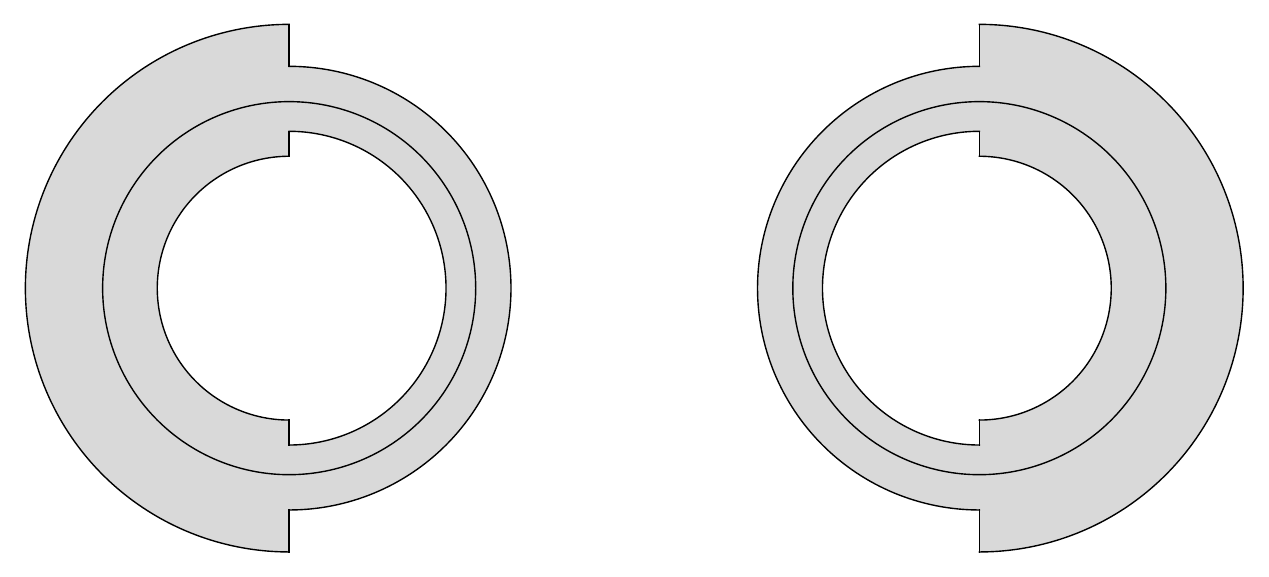}}
\unitlength1cm
\begin{picture}(0.01,0.01)
  \put(-10.65,-.35){$U_1$  }
  \put(-3.33,-.35){$U_2$  }
\end{picture}
\end{center}

\vspace{.7cm}

The reader might wonder if for a sequence $\{\Umum\}_{m=1}^\infty$ the convergence of meridians and their lengths together with that of the sequence $\{u_m\}_{m=1}^\infty$ is sufficient to ensure that $\{\Umum\}_{m=1}^\infty$ converges in the \Car topology. The example in Figure 3 shows that this is not the case, the basic reason being that knowing the meridians of a domain does not allow one to determine the domain itself. In both $U_1$ and $U_2$ the circle indicated is the unit circle and the domains are symmetric under $1/z$ so that by Theorem 1.7 this circle is the equator of the topological annulus concerned. However, it is clear that a sequence of pointed domains which alternated between these two could not converge in the \Car topology.

\section{Riemann Mappings}

In this section we prove a version of Theorem 1.2 for Riemann maps instead of covering maps. This is useful in situations where one wants to investigate properties of a family of functions where the functions are defined on different domains of the same connectivity. The usual thing to do is to normalize the domains to make them as similar as possible. However, given that even the normalized domains will likely be different, we need a notion of convergence of a sequence of functions on defined on varying domains. Of course, this is only likely to make sense if the domains themselves are also converging. 

\begin{definition}
Let $\{\Umum\}_{m=1}^\infty$ be a sequence of pointed domains which converges in the \Car topology to a pointed domain $\Uu$ with $\Uu \ne (\{u\},u)$. For each $m$ let $f_m$ be an analytic function (with respect to the spherical topology) defined on $U_m$ and let $f$ be an analytic function defined on $U$. We say that $f_m$ converges to $f$ {\rm uniformly on compact subsets of $U$} or simply {\rm  locally uniformly to $f$ on $U$} if, for every compact subset $K$ of $U$ and every $\epsilon >0$, there exists $m_0$ such that ${\mathrm d^\#}(f_m(z), f(z)) < \epsilon$ on $K$ for all $m \ge m_0$.
\end{definition}  

This is an adaptation to the spherical topology of the definition originally given in \cite{Ep}.
Note that, in view of condition ii) of \Car convergence, for any such $K$ $f_m$ will be defined on $K$ for all sufficiently large $m$ and so the definition is meaningful. Clearly if all the domains involved are the same, then we recover the standard definition of uniform convergence on compact subsets. This version of local uniform convergence is further related to the standard one in view of the following result whose proof is a straightforward application of Theorem 1.2 combined with ii) of \Car convergence. 

\begin{proposition} Let $\{\Umum\}_{m=1}^\infty$ be a sequence of pointed domains, which converges to $\Uu$ with $\Uu \ne (\{u\},u)$ and let $\pi_m$ and $\pi$ be the normalized covering maps from $\D$ to each $U_m$ and $U$ respectively. 
Let $f_m$ be defined on $U_m$ for each $m$ and $f$ be defined on $U$ and suppose 
$f_m$ converges uniformly to $f$ on compact subsets of $U$. Then the compositions 
$f_m \circ \pi_m$ converge locally uniformly on $\D$ to $f \circ \pi$.
\end{proposition}

Recall that there is a version of the Riemann mapping theorem for multiply connected domains which maps a given multiply connected domain into a domain of the same connectivity which is of some standard shape. There is some difference regarding the precise form of these standard domains: however, one of the most common is a round annulus from which a number of concentric circular slits have been removed such as can be found in Ahlfors' book \cite{Ahl}. From now on, we shall refer to such domains as \emph{standard domains} (where in the case $n=1$ the standard domain is the unit disc).

\begin{theorem}[\cite{Ahl} Page 255, Theorem 10] For an $n$-connected non-degenerate pointed domain $\Uu$ with $n > 1$, there is a conformal mapping $\phi(z)$ which maps $U$ to an annulus ${\mathrm A}(0,1,e^{\lambda^1})$ minus $n-2$ concentric arcs situated on the circles ${\mathrm C}(0,e^{\lambda^i})$, $i = 2, \ldots \ldots, n-1$. Furthermore, up to a choice of which complementary components of $U$ correspond to $\overline \D$ and 
$\cbar \setminus {{\mathrm D}(0,e^{\lambda^1})}$, the numbers $\lambda^i$, $i = 1, \ldots \ldots, n-1$ are uniquely determined as are the positions of the slits up to a rotation. If in addition we require that $\phi^\#(u) > 0$, the map $\phi$ is uniquely determined.
\end{theorem}

We remark that, despite the fact that the Riemann mapping does not in general extend beyond $U$, the construction Ahlfors gives shows how the correspondence between complementary components of $U$ and of the image domain can be done in a way which is both well-defined and natural.

We recall a well-known lemma concerning the behaviour of the hyperbolic metric near the boundary. A proof of the original version for the Euclidean metric can be found in \cite{CG} Page 13, Theorem 4.3 and it is the lower bound it gives on the hyperbolic metric which will be of particular importance for us. For a point $u \in U$, we shall denote the spherical distance to $\partial U$ by $\delta^\#_U (u)$ or just $\delta^\# (u)$ if once again the domain is clear from the context. 

\begin{lemma} Let $U \subset \cbar$ be a hyperbolic domain. Then there exists $C  > 0$ for which the hyperbolic metric $\rho(\cdot\,, \cdot)$ on $U$ satisfies
\[  \frac{C}{\delta^\#(z)\log(1/\delta^\#(z))} \dsharpz \le {\rm d}\rho (z) \le \frac{4}{\delta^\#(z)} \dsharpz, \qquad \mbox{as} \quad z \to \partial U.\]
\end{lemma}

The upper bound follows from the result in \cite{CG} combined with the facts that $\delta^\#(z)$ is less or equal to than the Euclidean distance to the boundary and that the spherical and Euclidean metrics are equivalent to within a factor of $2$ on the unit disc while the quantities $\dsharpz$, $\delta^\#(z)$ are invariant under the map $z \mapsto \tfrac{1}{z}$.

To obtain the lower bound, one lets $z_1$ be the closest point in $\cbar \setminus U$ to $z$ and chooses two other points $z_2$, $z_3$ in $\partial U$. These three points are then mapped using a M\"obius transformation to $0$, $1$ and $\infty$ respectively and one then obtains a lower bound on the hyperbolic metric for $\cbar \setminus \{0,1,\infty\}$ and applies the Schwarz lemma. It thus follows from Theorem 2.3.3 on page 34 of \cite{Bear} that these estimates are uniform with respect to the minimum separation in the spherical metric between $z_1$, $z_2$ and $z_3$.

Meridians are conformally invariant in the following sense. 

\begin{lemma} If $U$ is a hyperbolic domain and $\phi$ is a univalent function defined on $U$, then $\gamma$ is a meridian of $U$ if and only if $\phi(\gamma)$ is a meridian of $\phi(U)$. Furthermore, $\gamma$ is a principal meridian if and only if $\phi(\gamma)$ is.
\end{lemma}

\proof As before, we can assume that both $U$ and $\phi(U)$ are subdomains of $\C$.
$\gamma$ is a geodesic in $U$ if and only $\phi(\gamma)$ is a geodesic in $\phi(U)$. Also, two curves $\gamma_1$ and $\gamma_2$ are homologous in $U$ if and only if $\phi(\gamma_1)$ and $\phi(\gamma_2)$ are homologous in $\phi(U)$. The first part of the statement now follows from the conformal invariance of hyperbolic length.

For the second part, by invariance of homotopy or homology, if $\gamma$ is a simple closed curve in $U$, then $\gamma$ separates $\cbar \setminus U$ if and only if $\phi(\gamma)$ separates $\cbar \setminus \phi(U)$. It is then not too hard to see that by Theorem 1.3, we have that if $\gamma$ separates $\cbar \setminus U$ into two non-empty subsets $E$, $F$ and $\phi(\gamma)$ separates $\cbar \setminus \phi(U)$ into non-empty subsets $\tilde E$, $\tilde F$, then $E$ and $F$ are both disconnected if and only if $\tilde E$ and $\tilde F$ are. $\Box$

We will also need the following lemma on the conformal invariance of non-dgeneracy for finitely connected domains. 

\begin{lemma} Let $U$ be an $n$-connected domain with $n \ge 1$ and let $\phi$ be a univalent function defined on $U$. Then $U$ is non-degenerate if and only if $\phi(U)$ is. 
\end{lemma}

\vspace{-.31cm}

\proof For the case $n=1$, this is immediate from the Riemann mapping theorem in the simply connected case and the fact that $\C$ and $\D$ are not conformally equivalent. 

For $n \ge 2$, recall that a domain is degenerate if and only if we can find a 
curve in the domain which is homotopic to a puncture and contains curves of 
arbitrarily short hyperbolic length in its homotopy class. Since hyperbolic length and homotopy are both preserved by $\phi$, the result follows. $\Box$



Recall that a Riemann map to an $n$-connected slit domain as above with $n > 1$ is specified by $3n-5$ real numbers $\lambda^1, \lambda^2, \ldots, \lambda^{n-1}, \theta^1, \theta^2 \ldots, \theta^{2n-4}$ (we remark that Ahlfors considers the domains rather than the mappings, in which case, one can make an arbitrary rotation which allows one to eliminate one parameter in which case the domain is specified by $3n-6$ real numbers). Representing this list of numbers as a vector $\Lambda$, let us designate the pointed standard domain by $\ALa$ where the
inner radius is $1$, the outer radius $e^{\lambda^1}$ and the remaining $n-2$ complementary components are circular slits which are arcs of the circles ${\mathrm C}(0, e^{\lambda^j})$ which run from $e^{\lambda^j + i\theta^{2j - 3}}$ to $e^{\lambda^j + i\theta^{2j-2}}$. If $\Uu$ is mapped by the unique suitably normalized Riemann map $\phi$ to the pointed standard domain $\ALa$ where $\phi(u) = a$, $\phi'(u) >0$, we shall call $\ALa$ a \emph{standard domain for $U$}. 
This domain is unique up to assignment of which complementary components of $U$ get mapped to $\overline \D$ and $\cbar \setminus \overline{{\mathrm D}(0,e^{\lambda^1})}$.

Before stating the result, we remark that we only consider sequences of domains which have the same connectivity. To see why this is necessary, consider, for example, a pointed domain $(U,u)$ of low connectivity which is the limit of a sequence $\Umum$ where the domains $U_m$ have high connectivity which tends to infinity and where the diameters of the complementary components of $U_m$ all tend to zero. 

For $m$ large, at least one of the complementary components $L^i$ of $U$ is close (in the sense of the Euclidean or spherical distance distance between sets) to many complementary components of $U_m$.  However, this leads to two problems: firstly just which component of $\cbar \setminus U_m$ should one choose to correspond to a slit which is close to the corresponding slit for $L^i$ and secondly the fact that the components of $\cbar \setminus U_m$ could be very far apart relative to their size which could make the outer radius of $\ALm$ potentially very large (or even infinite) if one of these widely separated components corresponds to either of the components $\overline \D$ or $\cbar \setminus {\mathrm D}(0,e^{\lambda^1})$ of the complement of the standard domain $\ALm$. 

For a sequence of standard pointed domains $\{\ALam\}_{m=1}^\infty$, convergence in the \Car topology to another $n$-connected pointed domain $\ALa$ is precisely equivalent to the convergence of the points $a_m$ to $a$ and convergence in $\R^{3n -5}$ of the vectors $\Lambda_m$ to the corresponding vector $\Lambda$ for $\ALa$.
Finally, we remark that the behaviour and conformal invariance of the meridians and 
their lengths and the use of Theorem 1.10 are right at the heart of the proof of this result. Not surprisingly, Theorem 1.2 also plays a major role. 

\begin{theorem} Let $n \ge 1$, let $\{\Umum\}_{m=1}^\infty$ be a sequence of $n$-connected non-degenerate pointed domains, let $\Uu$ be an $n$-connected non-degenerate pointed domain and let $\ALa$ be a pointed standard domain for $\Uu$ where $a$ is the image of $U$ under the corresponding normalized Riemann map $\phi$ as in Theorem 3.1 (where we make any choice we wish regarding which components of $\cbar \setminus U$ correspond to $\overline \D$ and the unbounded complementary component of $A^\Lambda$). 

Then $\Umum$ converges to $\Uu$ if and only if we can label the components of the complements $\cbar \setminus U_m$ and choose corresponding normalized Riemann mappings $\phi_m$ to standard domains $\ALam$ so that these standard domains converge to $\ALa$ and the inverses $\psi_m$ of the maps $\phi_m$ converge locally uniformly on $\ALa$ to $\psi = \phi^{\circ -1}$, the inverse Riemann map for $\Uu$.

In addition, in the case of convergence, the 
Riemann maps $\phi_m$  converge locally uniformly on $\Uu$ to the Riemann map $\phi$ for $\Uu$. 
\end{theorem}

\proof The case $n=1$ is already proved in Theorem 1.2, so let us from now on assume that $n \ge 2$ and that the standard domains are then annuli from which (possibly) some slits have been removed. 

Suppose first that $\Umum$ converges to $\Uu$ and assume without loss of generality that $U \subset \C$.  The sequence $\{u_m\}$ of base points converges to $u$ and by discarding finitely many members if needed, we can assume that this sequence is bounded (in $\C$). Next, let $L^i$, $1 \le i \le n$ be the components of $\cbar \setminus U$ which correspond to our choice of standard domain (i.e. $L^1$ and $L^n$ correspond respectively to $\overline \D$ and $\cbar \setminus \overline {\rm D}(0, e^{\lambda^1})$). 

By the Hausdorff version of \Car convergence, any Hausdorff limit of the sets $\cbar \setminus U_m$ is contained in $\cbar \setminus U$. Using Lemma 2.1, we can label the components $\Kim$ of $\cbar \setminus U_m$ so that for each $1 \le i \le n$, and each component $L^i$ of $\cbar \setminus U$, any Hausdorff limit of the sets $\Kim$ is a subset of the component $L^i$ of $\cbar \setminus U$. 
 
We claim that the numbers 
$\lambda^1_m$ must be bounded above since otherwise, as each of the sets $\cbar \setminus \ALm$ has only $n$ components of which $n-2$ are slits, there would be a subsequence $m_k$ for which the standard annuli 
$A^{\Lambda}_{m_k}$ would contain round annuli whose moduli tended to infinity. By conformal invariance, we could say the same about the domains 
$U^{m_k}$ (where such thick annuli would separate the complements of these domains). The hyperbolic lengths of the equators of these annuli would then tend to $0$ and, by Theorem 1.6, the lengths of any meridians in the same homology classes as these equators would also tend to $0$. However, 
Corollary 1.1 tells us that the lengths $\leim$, $1 \le i \le E(n)$ of each $U_m$ are bounded below away from zero which then gives us a contradiction. 

By Montel's theorem and ii) of \Car convergence, the Riemann maps then give a normal family on any subdomain of $U$ which is compactly contained in $U$. A standard diagonal argument then shows that they must give a normal family on $U$ in the sense that any sequence taken from this family will have a subsequence which converges uniformly on compact subsets of $U$ in the sense of Definition 3.1. 

Now let $\gamma^i$, $1 \le i \le P(n)$ be the principal system of meridians which exists by Proposition 1.1 and using Lemma 3.2, we can consider the corresponding meridians $\tilde \gamma^i_m$, $1 \le i \le P(n)$ of $\ALm$. By relabelling if needed, we can say that the meridian $\tilde \gamma^1_m$ then separates $\overline \D$ from the rest of $\ALm$. Since the numbers $\lambda^1_m$ are uniformly bounded above, it follows that the spherical diameters of these meridians are bounded below. Additionally, by Theorem 1.10 and the conformal invariance of the hyperbolic metric, the lengths of these curves are uniformly bounded above. In view of our remarks after Lemma 3.1, we can use the estimates this result gives on the hyperbolic metric in a uniform fashion and since the improper integral

\vspace{-.4cm}
\[ \int_0^{\tfrac{1}{2}} { \frac{1}{x \log (1/ x)}{\rm d}x} \]

\vspace{.3cm}
diverges, we can find $\delta>0$ such that for every $m$ a spherical  $\delta$-neighbourhood of the meridian $\tilde \gamma^1_m$ is contained in $\ALm$. It then follows again by Theorem 1.10 and conformal invariance combined with the same estimates on the hyperbolic metric that we can make $\delta >0$ smaller if needed so that the spherical distance to the boundary $\delta^\#_{U_m}(a_m) \ge \delta$ for every $m$ and a spherical
$\delta$-neighbourhood of each of the principal meridians $\tilde \gamma^i_m$ of $\ALm$ will be contained in $\ALm$ for $1 \le i \le P(n)$ and every $m$.  In particular this means that the complementary components of each $\ALm$ are at least $2\delta$ away from each other. 

By the Koebe one-quarter theorem, the absolute values of the derivatives $\phi_m'(u_m)$ are uniformly bounded above and below away from $0$, whence all limit functions for the sequence $\{\phi_m\}_{m=1}^\infty$ must be non-constant and in fact univalent in view of Hurwitz's theorem. We next want to show that the standard domains $\ALam$ converge in the \Car topology and we will do this by appealing to Theorem 1.2.

Suppose that we can find a subsequence $m_k$ for which $\phi_{m_k}$ converges on $(U,u)$ to some univalent limit function $\phi$. Recall the normalized covering maps $\pi_m : \D \to U_m$ of Theorem 1.2 which by this result converge to the normalized covering map $\pi: \D \to U$. 

If we now set $\chi_{m_k} = \phi_{m_k} \circ \,\pi_{m_k}$, then $\chi_{m_k}$ is the unique normalized covering map for the standard pointed domain $(A^{\Lambda_{m_k}},a_{m_k})$. By Proposition 3.1, the functions $\chi_{m_k}$ then converge on compact subsets of $\D$ to $\phi \circ \pi$. Since $\phi$ is univalent and $\pi$ is a covering map, $\phi \circ \pi$ is itself a covering map which must in fact be $\chi$, the normalized covering map from $\D$ to $\phi(U)$. 

By Theorem 1.2, the domains $(A^{\Lambda_{m_k}},a_{m_k})$ then converge to a limit domain $(A',a')$ where $A' = \phi (U)$, and since $\delta(a_m) \ge \delta$ and the complementary components of each $\ALm$ are at least $2\delta$ away from each other for every $m$, $(A',a') \ne \{a'\}$ and must be $n$-connected. Also, from the Hausdorff version of \Car convergence, it follows that $(A',a')$ must be a standard pointed domain. Finally, as $U$ is non-degenerate and $\phi$ is univalent, it follows again by Lemma 3.3 that $A'$ is also non-degenerate. 

Now $\phi$ is univalent on $U$ and clearly $\phi^\#(u) >0$, so $\phi$ is the normalized Riemann map from $\Uu$ to $(A',a')$. By Theorem 3.1, $A'$ is conformally equivalent to $A^{\Lambda}$ and in order to show these two domains are equal we just need to show that 
$\phi$ preserves the labelling of the components of $\cbar \setminus U$.

Let $\gamma$ be a simple closed curve around the complementary component $L^1$ of $U$ which does not encircle the other complementary components of $U$ and which exists in view of Theorem 1.3. By our labelling of the complementary components of the domains $U_m$ and ii) of \Car convergence, $\gamma$ separates $K^1_{m_k}$ from the other components of $\cbar \setminus U_{m_k}$ for $k$ large enough. From this it is not too hard to see that, for $k$ large enough, $\phi_{m_k}(\gamma)$ is then a simple closed curve which separates $\overline \D$ from the other components of $\cbar \setminus A^{\Lambda_{m_k}}$ and thus encloses $\overline \D$. If we then let $z$ be any point of $\overline \D$, then by the local uniform convergence of $\phi_{m_k}$ to $\phi$ on $U$, $n(\phi(\gamma),z) = n(\phi_{m_k}(\gamma),z) = \pm 1$ for $k$ large enough whence $\phi(\gamma)$ also encloses $\overline \D$. It also follows from Lemma 2.1 and the convergence of the pointed domains $(A^{\Lambda_{m_k}},a_{m_k})$ to $(A',a')$ 
that $\phi(\gamma)$ does not enclose any of the other components of $\cbar \setminus A'$. 

$L^1$ thus corresponds under $\phi$ to the complementary component $\overline \D$ of $A'$ (and also $A^{\Lambda}$) and a similar argument shows that $L^n$ corresponds to the unbounded complementary component of $A'$. By the uniqueness part of Theorem 3.1, we must then have that $(A',a') = \ALa$. It then follows easily that $\ALam$ converges to $\ALa$ and that the mappings $\phi_m$ converge on compact subsets of $U$ to $\phi$. 

We still need to show the inverses $\psi_m$ converge. Since the domains $\Umum$ converge to another pointed $n$-connected domain none of whose complementary components is a point, by Lemma 2.1 the spherical diameters of the complements $\cbar \setminus U_m$ are bounded below and the usual argument of post-composing with a bi-equicontinuous family of M\"obius transformations and applying Montel's theorem shows that the functions $\psi_m$ give a family which is normal on $A^\Lambda$ in the sense given earlier.

Applying the Koebe one-quarter theorem and Hurwitz's theorem as before shows that all limit functions must be non-constant and univalent. If we then have a sequence $\psi_{m_k}$ which converges uniformly on compact subsets of $A^\Lambda$ to a limit function $\tilde \psi$, then, by Proposition 3.1 again, $\psi_{m_k} \circ \chi_{m_k}$ converges uniformly on compact subsets of $\D$ to $\tilde \psi \circ \chi$. Using Rouch\'e's theorem and local compactness as in \cite{Ep} then shows that $\tilde \psi (A^\Lambda) = U$ with $\tilde \psi (a) = u$ and using ii) of \Car convergence, it follows easily that $\phi_{m_k} \circ \psi_{m_k}$ converges uniformly on compact subsets of $A^\Lambda$ to $\phi \circ \tilde \psi$ whence $\tilde \psi = \phi^{\circ -1}$. With this the proof of the first direction is finished. 

For the other direction, suppose now that the standard pointed domains $\ALam$ converge to $\ALa$, which is an $n$-connected non-degenerate standard domain and that the corresponding inverse Riemann maps $\psi_m$ converge to $\psi$. For each $m$ let $\chi_m$ be the normalized covering map from $\D$ to the standard domain $\ALm$ and let $\chi$ be the corresponding covering map for $\AL$ so that $\chi_m$ converges uniformly on compact subsets of $\D$ to $\chi$ by Theorem 1.2. Again by Proposition 3.1, $\pi_m = \chi_m \circ \psi_m$ will converge 
locally uniformly to $\chi \circ \psi = \pi$ on $\D$. By Theorem 1.2, it then follows that 
$\Umum$ converges to $\Uu$. On the other hand, as $\Uu$ is a \Car limit of pointed domains of connectivity $n$, $U$ has connectivity $\le n$ and since $A^\Lambda$ is $n$-connected and $\psi$ is univalent, it follows from Theorem 3.1 that $U$ must be $n$-connected.  Finally, as $A^\Lambda$ is non-degenerate, it follows from Lemma 3.3 that $U$ must be non-degenerate. $\Box$


\end{document}